\input ./style/arxiv-general.cfg
\documentclass[aop,nameyear,MSNbibl,seceqn,dvips]{arximspdf}
\makeatletter
   \@ifpackageloaded{graphicx}{}{\usepackage{graphicx}}
\makeatother

\doi{10.1214/14-AOP961} 
\volume{43}
\issue{6}
\pubyear{2015}
\firstpage{3337}
\lastpage{3358}
\docsubty{FLA}

\makeatletter
\newcommand{\aaa}{\operatorname{\mathtt{cogr}}}
\renewcommand{\mid}{\vert}
\newproclaim{remark}[theorem]{Remark}
\newtheorem{theorem}{Theorem}[section]
\newtheorem{proposition}[theorem]{Proposition}
\newtheorem{lemma}[theorem]{Lemma}
\makeatother

\begin{document}
\begin{frontmatter}

\title{Cycle density in infinite Ramanujan graphs\thanksref{T1}}
\runtitle{Cycle density in infinite Ramanujan graphs}

\begin{aug}
\author[A]{\fnms{Russell}~\snm{Lyons}\corref{}\ead[label=e1]{rdlyons@indiana.edu}\ead[label=u1,url]{http://pages.iu.edu/\textasciitilde rdlyons/}}
\and
\author[B]{\fnms{Yuval}~\snm{Peres}\ead[label=e2]{peres@microsoft.com}\ead[label=u2,url]{http://research.microsoft.com/en-us/um/people/peres/}}
\runauthor{R. Lyons and Y. Peres}
\affiliation{Indiana University and Microsoft Corporation}
\address[A]{Department of Mathematics\\
Indiana University\\
831 E 3rd St\\
Bloomington, Indiana 47405-7106\\
USA\\
\printead{e1}\\
\printead{u1}} 
\address[B]{Microsoft Corporation\\
One Microsoft Way\\
Redmond, Washington 98052-6399\\
USA\\
\printead{e2}\\
\printead{u2}}
\end{aug}

\thankstext{T1}{Supported in part by Microsoft Research and
NSF Grant DMS-10-07244.}

\received{\smonth{10} \syear{2013}}
\revised{\smonth{8} \syear{2014}}

%
\begin{abstract}
We introduce a technique using nonbacktracking random walk for
estimating the spectral radius of simple random walk.
This technique relates the density of nontrivial cycles in simple random
walk to that in nonbacktracking random walk. We apply this to
infinite Ramanujan graphs, which are regular graphs whose spectral
radius equals that of the tree of the same degree.
Kesten showed that the only infinite Ramanujan graphs that are Cayley
graphs are trees.
This result was extended to unimodular random rooted regular graphs by
Ab\'ert, Glasner and Vir\'ag.
We show that an analogous result holds for all regular graphs: the
frequency of times spent by simple random walk in a nontrivial cycle is
a.s. 0 on every infinite Ramanujan graph.
We also give quantitative versions of that result, which we apply to answer
another question of
Ab\'ert, Glasner and Vir\'ag, showing that on an infinite Ramanujan graph,
the probability that simple random walk encounters a short cycle tends
to 0
a.s. as the time tends to infinity.
\end{abstract}

%
\begin{keyword}[class=AMS]
\kwd[Primary ]{60G50}
\kwd{82C41}
\kwd[; secondary ]{05C81}.
\end{keyword}
\begin{keyword}
\kwd{Nonbacktracking random walks}
\kwd{spectral radius}
\kwd{regular graphs}.
\end{keyword}
\end{frontmatter}

\section{Introduction}\label{s.intro}

A path in a multigraph is called \textit{nonbacktracking} if no edge is
immediately followed by its reversal.
Note that a loop is its own reversal.
Nonbacktracking random walks are almost as natural as ordinary
random walks, though more difficult to analyze in most situations.
Moreover, they can be more useful than ordinary random walks when random
walks are used to search for something, as they explore more quickly, not
wasting time immediately backtracking; see \citet{MR2348845}.
Our aim, however, is to use
them to analyze the spectral radius of ordinary random walks on regular
graphs.

The \textit{spectral radius} of a (connected, locally finite)
multigraph $G$
is defined to be $\rho(G) := \limsup_{n \to\infty} p_n(o, o)^{1/n}$
for a
vertex $o \in G$, where $p_n(x, y)$ is the $n$-step transition probability
for simple random walk on $G$ from $x$ to $y$. It is well known that
$\rho(G)$ does not depend on the choice of $o$.

If $G = \mathbb{T}_d$ is a regular tree of degree $d$, then $\rho(G) =
2\sqrt{d-1}/d$. Regular trees are Cayley graphs of groups. In general,
when $G$ is a Cayley graph of a group, \citet{MR22253}
proved that $\rho(G) > \rho(\mathbb{T}_d)$ when $G$ has degree $d$
and $G \ne
\mathbb{T}_d$. \citet{MR222911} also proved that for Cayley graphs,
$\rho(G) =1$ iff $G$ is amenable.

If $G$ is a $d$-regular multigraph, then its universal cover is
$\mathbb{T}_d$, whence $\rho(\mathbb{T}_d) \le\rho(G) \le1$.
Using the method of proof
due to \citet{MR536645}, various researchers related $1 - \rho(G)$ to the
expansion (or isoperimetric) constant of infinite graphs $G$, showing that
again, $G$ is amenable iff $\rho(G) = 1$; see \citet{MR85m58185}, \citet{MR88h58118}, \citet{MR87e60124}, \citet{MR973878}, \citet{MR89g60216},
\citet{MR89a05103} and \citet{MR94i31012}.

It appears considerably more difficult to understand the other inequality
for $\rho(G)$: when is $\rho(G) = \rho(\mathbb{T}_d)$?
This question will be our focus.

For finite graphs, the spectral radius is 1. Of interest instead is the
second largest eigenvalue, $\lambda_2$, of the transition matrix.
An inequality of Alon and Boppana [see \citet{MR88e05077} and \citet{MR1124768}] says that if $\langle{G_n; n \ge1}\rangle$
is a
family of
$d$-regular graphs whose size tends to infinity, then $\liminf_{n
\to\infty} \lambda_2(G_n) \ge\rho(\mathbb{T}_d)$.
Regular graphs $G$ such that all eigenvalues have absolute value either 1
or at most $\rho(\mathbb{T}_d)$ were baptized \textit{Ramanujan graphs} by
\citet{MR963118}, who, with \citet{MR939574},
were the first to exhibit explicit such families. Moreover,
their examples had better expansion properties than the random graphs that
had been constructed earlier. See \citet{MR1966527}
and \citet{MR2331350} for
surveys of finite Ramanujan graphs.

\citet{AGVmeasurable} studied the density of short
cycles in Ramanujan
graphs. One of their tools was graph limits, which led them to define and
study \textit{infinite Ramanujan graphs}, which are
$d$-regular infinite
graphs whose spectral radius equals $\rho(\mathbb{T}_d)$.
Now limits of finite graphs, taken in the appropriate sense, are
probability measures on rooted graphs; the probability measures that arise
have a property called unimodularity.
Theorem~5 of \citet{AGVmeasurable} shows that every
unimodular
random rooted infinite regular graph that is a.s. Ramanujan is a.s. a tree.
Unimodularity is a kind of stochastic homogeneity that, among other things,
ensures that simple random walk visits short
cycles with positive frequency when they exist.

\citet{AGVmeasurable} asked whether the hypothesis
of unimodularity could
be weakened to something called stationarity. We answer this affirmatively
in a very strong sense, using no extra hypotheses on the graph and
including cycles of all lengths at once.
To state our result,
call a cycle \textit{nontrivial} if it is not purely a
backtracking cycle,
that is, if when backtracks are erased iteratively from the cycle, some edge
remains.
For example, a single loop is a nontrivial 1-edge cycle, but a~loop
followed by the same loop is a trivial 2-edge cycle.
Let $X = \langle{X_n; n \ge 1}\rangle$ be simple random walk on $G$,
where $X_n$
are directed edges and the tail of $X_1$ is any fixed vertex.
Call $n$ a \textit{nontrivial cycle time} of $X$ if there
exist $1 \le s \le n \le t$ such that $(X_s, X_{s+1}, \ldots, X_t)$ is a nontrivial cycle.

\begin{theorem}\label{t.main}
If $G$ is an infinite Ramanujan graph of degree at least 3, then a.s. the
density of nontrivial cycle times of $X$ in $[1, n]$ tends to 0 as $n
\to\infty$.
\end{theorem}

Now fix $L \ge 1$. Let $q_n$ be
the probability that simple random walk at time $n$ lies on a nontrivial
cycle of length at most $L$. Then the preceding theorem implies that
$\liminf_{n \to\infty} q_n = 0$.
In their Problem~10, \citet{AGVmeasurable} ask
whether $\lim_{n \to\infty} q_n = 0$.
We answer it affirmatively.

\begin{theorem}\label{t.prob10}
Let $G$ be an infinite Ramanujan graph and $L \ge 1$.
Then $\lim_{n \to\infty} q_n = 0$.
\end{theorem}

In broad outline,
our technique to prove these results is the following:
First, we prove that when simple
random walk on $G$ has many nontrivial cycle times, then so does
nonbacktracking random walk.
Second, we deduce that under these circumstances, we may transform
nonbacktracking paths to nonbacktracking cycles with controlled length
and find that there are many nonbacktracking cycles.
The exponential growth rate of the number of nonbacktracking cycles is
called the cogrowth of $G$.
Finally, we use the cogrowth formula relating cogrowth to spectral radius
to conclude that $G$ is not Ramanujan.

Thus, of central importance to us is the notion of cogrowth. We state the
essentials here.
Let the number of nonbacktracking cycles of length $n$ starting from some
fixed $o \in\mathtt{V}(G)$ be $b_n(o)$.
Let
\[
\aaa (G) := \limsup_{n \to\infty} b_n(o)^{1/n}
\]
be the exponential
growth rate of the number of nonbacktracking cycles containing~$o$.
This number is called the \textit{cogrowth} of $G$.
The reason for this name is that if we consider a universal covering map
$\varphi\dvtx T \to G$, then the cogrowth of $G$ equals the exponential
growth rate of $\varphi^{-1}(o)$ inside $T$ since $\varphi$
induces a bijection between simple paths in $T$ and
nonbacktracking paths in $G$.
By using this covering map, one can see that $\aaa (G)$ does
not depend on
$o$.
Note, too, that if $\mathcal{P}$ is a finite path in $G$ that lifts to
a path in
$T$ from vertex $x$ to vertex $y$, then erasing backtracks from~$\mathcal{P}$
iteratively yields $\varphi[\mathcal{P}']$, where $\mathcal{P}'$ is
the shortest path in
$T$ from $x$ to $y$.

Let $G$ be a connected graph.
It is not hard to check the following:
If $G$ has no simple nonloop cycle and at most one loop,
then $\aaa (G) = 0$.
If $G$ has
one simple cycle and no loop or no simple nonloop cycle and
two loops, then $\aaa (G) = 1$.
In all other cases, that is, when the fundamental group of $G$ is not
virtually abelian, $\aaa (G) > 1$.

The central result about cogrowth is the following formula (\ref{e.cogrowth}),
due to \citet{MR599539} for Cayley graphs and extended by
\citet{MR1120509} to all regular graphs.

\begin{theorem}[(Cogrowth formula)]\label{t.cogrowth}
If  $G$ is a $d$-regular connected multigraph, then
%
\begin{equation}
\label{e.cogr-cond} \aaa (G) > \sqrt{d-1} \quad\mbox{iff}\quad \rho(G) >
{2\sqrt{d-1} \over d},
\end{equation}
in which case
%
\begin{equation}
\label{e.cogrowth} d \rho(G) = {d-1 \over\aaa (G)} + \aaa (G).
\end{equation}
If (\ref{e.cogr-cond}) fails, then
$\rho(G) = 2\sqrt{d-1}/d$
and $\aaa (G) \le\sqrt{d-1}$.
\end{theorem}

See \citet{LPbook}, Section~6.3, for a proof.

Our use of Theorem~\ref{t.cogrowth} will be mainly via (\ref{e.cogr-cond}), rather
than (\ref{e.cogrowth}).
In order to use (\ref{e.cogr-cond}), we shall prove the following
result on
density of nontrivial cycle times:

\begin{theorem}\label{t.main2}
Suppose that $G$ is a graph all of whose degrees are at least~3.
If with positive probability the limsup
density of nontrivial cycle times of simple random walk in $[1, n]$ is
positive as $n \to\infty$, then the same holds for nonbacktracking random
walk.
\end{theorem}

Here, nonbacktracking random walk is the random walk that at every
time~$n$, chooses
uniformly among all possible edges that are not the reversal of the $n$th
edge.
In terms of the universal cover $\varphi\dvtx T \to G$, if simple random
walk $X = \langle{X_n; n \ge1}\rangle$ is lifted to a random walk,
call it
$\widehat{X}$, on
$T$, then $\widehat{X}$ is simple random walk on $T$.
Backtracking on $G$ is the same as backtracking on $T$.
Since all degrees of $T$
are at least 3, $\widehat{X}$ is transient and so there is a unique
simple path
$\mathcal{P}$ in~$T$ with the same starting point as $\widehat{X}$
and having
infinite intersection with $\widehat{X}$. The law of
$\varphi[\mathcal{P}]$ is that of nonbacktracking random walk on $G$.

As it may be of separate interest, we note in passing the following basic
elementary bound on the number of nonbacktracking cycles.
Write $S(x) := \{ n; b_n(x) \ne0\}$.
If a nonbacktracking cycle is a loop or
has the property that its last edge is different from the
reverse of its first edge, then call the cycle \textit{fully
nonbacktracking} (usually called ``cyclically reduced'' in the case of a
Cayley graph).
Let the number of fully nonbacktracking
cycles of length $n$ starting from $x$ be $b_n^*(x)$.

\begin{proposition}\label{p.cogrlimit}
Let $G$ be a graph with $\aaa (G) \ge1$.
For each $x \in\mathtt{V}(G)$, we have that
$\lim_{S(x) \ni n \to\infty} b_n(x)^{1/n}$ exists and
there is
a constant $c_x$ such that $b_n(x) \le c_x \aaa (G)^n$ for all
$n \ge1$.
Furthermore, if $x$ belongs to a simple cycle of length $L$, then $c_x
\le
2 + 2L\aaa (G)^{L-2}$.
If $G$ is $d$-regular, then $G$ is Ramanujan iff for all vertices $x$ and
all $n \ge1$, we have $b_n^*(x) \le2 (d-1)^{n/2}$.
\end{proposition}


We shall illustrate our technique first
by giving a short proof of Kesten's theorem (extended to transitive
multigraphs).
We then prove a version of Theorems \ref{t.main} and~\ref{t.main2} with
a stronger hypothesis on the
density of nontrivial cycle times, a hypothesis that holds for stationary
random rooted graphs, for example.
The proof of the full Theorems \ref{t.main} and \ref{t.main2}
requires a large number of technical lemmas,
which makes the basic idea harder to see.
The final section proves Theorem~\ref{t.prob10}.

All our graphs are undirected connected infinite
multigraphs. However, each edge comes with
two orientations, except loops, which come with only one orientation.
An edge $e$ is oriented from its tail $e^{-}$ to its head $e^{+}$.
These endpoints are the same when $e$ is a loop.
A vertex may have many loops and two vertices may be joined by many edges.
If $e$ is an oriented edge, then its reversal is the same unoriented
edge with the opposite orientation, denoted $- e$. This is equal to
$e$ iff $e$ is a loop.

We shall have no need of unimodularity or stationarity, so we do not define
those terms.

\section{Kesten's theorem}\label{s.kesten}

It is easiest to understand the basic ideas behind our proofs in the case
of transitive multigraphs.
\citet{MR22253} proved the following result and
various extensions for Cayley graphs.

\begin{theorem}\label{t.kestenramanujan}
If $d \ge3$ and $G$ is a $d$-regular transitive multigraph that is not a
tree, then $\rho(G) > \rho(\Bbb T_d)$.
\end{theorem}


\begin{pf}
Let $L$ be the length of the shortest cycle in $G$ (which is 1 if
there is a loop).
Consider a nonbacktracking random walk $\langle{Y_n; n \ge1}\rangle
$, where each
edge $Y_{n+1}$ is chosen uniformly among the edges incident to the head
$Y_n^{+}$ of $Y_n$, other than the reversal of $Y_n$.
We are going to handle loops differently than other cycles, so it will be
convenient to let
\[
L' := \cases{L, & \quad \mbox{if }$L > 1$,
\cr
3, & \quad \mbox{if
}$L = 1$.}
\]
Let $A_n$ be the event that
$Y_{n+1}, \ldots, Y_{n+L'}$ is a nonbacktracking cycle.
Write $b := d-1$.
For $n \ge 1$,
\[
\mathbf{P}(A_n \mid Y_1, \ldots, Y_n) \ge
{1 \over d b^{L'-1}}
\]
since if $L > 1$, then there is a way to traverse a simple cycle starting
at $Y_n^{+}$ and not using the reversal of $Y_n$,
while if $L = 1$, then the walk can first choose an edge other than the
reversal of $Y_n$, then traverse a loop, and then return by the
reversal of
$Y_{n+1}$.
Let $Z_k := \mathbf{1}_{A_{k L'}} - \mathbf{P}(A_{k L'} \mid Y_1,
\ldots, Y_{k L'})$.
Then $\langle{Z_k; k \ge1}\rangle$ are uncorrelated, whence by the
Strong Law of
Large Numbers for uncorrelated random variables, we have
\[
\lim_{n \to\infty} {1 \over n} \sum
_{k=0}^{n-1} Z_k = 0 \qquad\mbox{a.s.},
\]
which implies that
\[
\liminf_{n \to\infty} {1 \over n} \sum
_{k=0}^{n-1} \mathbf{1}_{A_{k L'}} \ge
{1 \over d b^{L'-1}} \qquad\mbox{a.s.}
\]
Therefore, if we choose $\varepsilon< 1/  (d b^{L'-1} )$, then in
$n L'$
steps, at least $\varepsilon n$ events $A_{k L'}$ will occur for $0 \le k <
n$ with probability tending to 1 as $n \to\infty$.

Consider the following transformation of a path $\mathcal{P}= (Y_1,
\ldots,Y_{nL'})$ to a ``reduced'' path $\mathcal{P}'$:
For each $k$ such that $A_{kL'}$ occurs,
remove the edges $Y_{k+1}, \ldots,  Y_{k+L'}$.
Next, combine $\mathcal{P}$ and $\mathcal{P}'$ to form a
nonbacktracking cycle $\mathcal{P}''$ by
appending to $\mathcal{P}$ a nonbacktracking cycle of length
$L'$ that does not begin with the reversal of $Y_{nL'}$,
and then by returning to the tail of $Y_1$ by $\mathcal{P}'$
in reverse order. 
Note that the map $\mathcal{P}\mapsto\mathcal{P}''$ is 1--1.

When at least $\varepsilon n$ events $A_{k L'}$ occur,
the length of $\mathcal{P}''$ is at most $(2n + 1 - \varepsilon n)L'$.
The number of nonbacktracking paths $Y_1, \ldots, Y_n$ equals $d b^{n-1}$,
whence $\sum_{k \le(2n + 1 - \varepsilon n)L'} b_k(G) \ge d b^{nL'-1}/2$
for\vspace*{1pt} large $n$.
This gives that\break  $\aaa (G) > \sqrt b$, which implies the result by
Theorem~\ref{t.cogrowth}.
\end{pf}

An alternative way of handling loops in the above proof is to use
the following:
Consider a random walk on a graph with spectral radius $\rho$.
Suppose that we introduce a delay so that each step goes nowhere with
probability $p_{\mathrm{delay}}$, and otherwise chooses a neighbor with the same
distribution as before. Then the new spectral radius equals $p_{\mathrm{delay}}+ (1-p_{\mathrm{delay}}) \rho$.
Hence, if there is a loop at each vertex and $G$ is $d$-regular, then
$\rho(G) \ge1/(d-1) + (d-2)\rho(\mathbb{T}_d)/(d-1) > \rho(\mathbb{T}_d)$.

For a simple extension of this proof, let $G$ be a $d$-regular multigraph.
Suppose that there are some $L, M < \infty$ such that for
every vertex $x \in\mathtt{V}(G)$, there is a simple cycle of length
at most
$L$ that is at distance at most $M$ from $x$. Then $\rho(G) > \rho
(\mathbb{T}_d)$.
Theorem~3 of \citet{AGVmeasurable} gives a
quantitative strengthening of
this result.

\section{Expected frequency}\label{s.expfreq}

In the case of transitive multigraphs that are not trees,
it is clear that simple random walk
a.s. has many nontrivial cycle times.
The most difficult part of our extension to general regular graphs is to
show how this property is inherited by nonbacktracking random walk.
This actually does not depend on regularity and is an interesting fact in
itself.

Before we prove the general case, which has many complications, it may be
helpful to the reader to see how to prove Theorem~\ref{t.main} with a stronger
assumption on the density of nontrivial cycle times.


Recall that a cycle is
\textit{nontrivial} if it is not purely a backtracking cycle,
that is, when backtracks are erased iteratively from the cycle, some
edge remains.
We call such cycles \textit{NT-cycles}.


\begin{theorem}\label{t.AGV}
Suppose that $G$ is a graph all of whose degrees lie in some
interval $[3, D]$.
If the limsup expected frequency that simple random walk traverses
some nontrivial cycle of length at
most $L$ is positive, then the same is
true for nonbacktracking random walk. Hence, if $G$ is also
$d$-regular, then $\rho(G) > 2\sqrt{d-1}/d$.
\end{theorem}

\begin{pf}
We may assume that
simple cycles of length exactly $L$ are traversed with positive
expected frequency.
Let $X = \langle{X_n; n \ge1}\rangle$ be simple random walk on $G$
and $\widehat{X}=
\langle{\widehat{X}_n}\rangle$ be its lift to the universal cover
$T$ of $G$.

Now consider $X$.
It contains purely backtracking excursions that are erased when we
iteratively erase all backtracking.
Let the lengths of the successive excursions be $M_1, M_2, \ldots\,$, where
$M_i \ge0$. Define
%
%
\begin{equation}\label{e.defPhi}
\Phi(n) := n + \sum_{k = 1}^{n}
M_k.
\end{equation}
Then the edges that remain
after erasing all backtracking are $\langle{X_{\Phi(n)}; n \ge
1}\rangle$.
If we write $Y_n := X_{\Phi(n)}$, then
$Y := \langle{Y_n}\rangle=: \mathtt{NB}(X)$ is the nonbacktracking
path created from $X$.
Let $\operatorname{im}\Phi$ be the image of $\Phi$.
Thus, $t \in\operatorname{im}\Phi$ iff the edge $X_t$ is not erased from $X$
when erasing all backtracking.

Consider a time $t$ such that $X_t$ completes a
traversal of a simple cycle of length~$L$.
Because all degrees of $T$ are at least 3, the probability (given the past)
that $\widehat{X}$ will never cross the edge $-\widehat{X}_t$ after
time $t$ is at least $1/2$.
In such a case, the cycle just traversed will not be erased (even in part)
by the future.
However, erasing backtracks from $X_1, \ldots, X_t$ may erase
(at least in part) this cycle.

Let $\mathtt{Trav}(Y)$ be the set of times $n$ such that $Y_n$
completes a traversal of a cycle
and let $\mathtt{Trav}(X)$ be the set of times $t$ that $X_t$
completes a traversal of a simple cycle of length $L$.

We divide the rest of the proof into two cases, depending on whether $L >
1$ or not.

First, suppose that $L > 1$.
Define a map $\psi\dvtx\mathbb{Z}^+ \to\mathtt{Trav}(Y) \cup\{
\infty\}$ as follows:
\[
\psi(t) := \cases{\Phi^{-1}(t)+L,  &  \quad \mbox{if }$t \in\operatorname{im}\Phi \cap
\mathtt{Trav}(X)$\mbox{ and }$\Phi^{-1}(t)+L \in\mathtt{Trav}(Y)$,
\cr
\infty, &\quad  \mbox{otherwise}.}
\]
For $t \in\mathtt{Trav}(X)$, the probability (given the past) that
the steps
$X_{t+1},
\ldots,\break  X_{t+L}$ traverse the same cycle $X_{t-L+1}, X_{t-L+2},
\ldots,
X_{t}$ of length\vspace*{1pt} $L$ (in only $L$ steps
and in the same direction) and then (on the tree)
$\widehat{X}_{t+L+1}, \widehat{X}_{t+L+2}, \ldots$ never
crosses the edge $-\widehat{X}_t$ is at least $1/(2D^L)$;
similarly, for traversing the cycle in the opposite direction.
In at least one of these two cases, some part of
the cycle $X_{t+1}, \ldots, X_{t+L}$ will
be left after erasing \textit{all} backtracks in $X$, in which case
$\psi(t)
\in\mathtt{Trav}(Y)$.
Therefore, $\mathbf{P} [\psi(t) \in\mathtt{Trav}(Y) | t
\in\mathtt{Trav}(X) ] \ge
1/(2D^L)$, that is,
$\mathbf{P} [\psi(t) \in\mathtt{Trav}(Y) ] \ge\mathbf
{P}[t \in\mathtt{Trav}(X)]/(2D^L)$.
Hence,
\[
\sum_{s \le t} \mathbf{P} \bigl[\psi(s) \in
\mathtt{Trav}(Y) \bigr] \ge \sum_{s \le t} \mathbf{P}
\bigl[s \in\mathtt{Trav}(X)\bigr]/\bigl(2D^L\bigr).
\]
Note that $\psi(t_1) = \psi(t_2) \in\mathtt{Trav}(Y)$ implies that
$t_1 = t_2$.
Since $\psi(s) \le s+L$, it follows that
\[
\sum_{k \le t+L} \mathbf{1}_{[{k \in\mathtt{Trav}(Y)}]} \ge \sum
_{s \le t} \mathbf{1}_{[\psi(s) \in\mathtt{Trav}(Y)]},
\]
whence
\begin{eqnarray*}
&&\limsup_{n \to\infty} n^{-1} \sum
_{k \le n} \mathbf{P}\bigl[k \in \mathtt{Trav}(Y)\bigr]
\\
&&\qquad\ge \limsup_{t \to\infty} t^{-1} \sum
_{s \le t} \mathbf{P}\bigl[s \in \mathtt{Trav}(X)\bigr]/
\bigl(2D^L\bigr)
\\
&&\qquad= \limsup_{t \to\infty} t^{-1} \mathbf{E} \biggl[
\sum_{s \le t} \mathbf{1}_{[s \in\mathtt{Trav}(X)]} \biggr] \Big/
\bigl(2D^L\bigr) > 0
\end{eqnarray*}
by assumption.
Now the method of proof of Theorem~\ref{t.kestenramanujan} applies when
$G$ is regular.

Finally, suppose that $L = 1$.
This means that $t \in\mathtt{Trav}(X)$ iff $X_t$ is a loop, and
similarly for
$\mathtt{Trav}(Y)$.
Define a map $\psi\dvtx\mathbb{Z}^+ \to\mathbb{Z}^+ \cup\{\infty
\}$ as follows:
\[
\psi(t) := \cases{\Phi^{-1}(t), & \quad \mbox{if }$t \in\operatorname{im}\Phi \cap
\mathtt{Trav}(X)$,
\cr
\Phi^{-1}(t+1),  &\quad \mbox{if }$t \in\mathtt{Trav}(X)
\setminus\operatorname{im}\Phi$ and $t+1 \in\operatorname{im}\Phi\cap\mathtt{Trav}(X)$,$\!\!$
\cr
\infty, & \quad $\mbox{otherwise}$.}
\]
Consider $t \in\mathtt{Trav}(X)$.
If erasing backtracks from $X_1, \ldots, X_t$ does not erase
$X_t$, then the probability (given the past) that
(on the tree) $\widehat{X}_{t+1}, \widehat{X}_{t+2}, \ldots$ never
crosses the edge $-\widehat{X}_t$ is at least $1/2$, in which case $t
\in\operatorname{im}\Phi$.
On the other hand,
if erasing backtracks from $X_1, \ldots, X_t$ \textit{does} erase
$X_t$,
then the probability (given the past) that $X_{t+1}$ is a loop
and (on the tree)
$\widehat{X}_{t+2}, \widehat{X}_{t+3}, \ldots$ never
crosses the edge $-\widehat{X}_{t+1}$ is at least $1/(2D)$, in which case
$t \notin\operatorname{im}\Phi$ and $t+1 \in\operatorname{im}\Phi\cap\mathtt
{Trav}(X)$.
In each of these two cases, $\psi(t) \in\mathtt{Trav}(Y)$.
Therefore, $\mathbf{P} [\psi(t) \in\mathtt{Trav}(Y) | t
\in\mathtt{Trav}(X) ] \ge1/(2D)$,
that is,\break  $\mathbf{P} [\psi(t) \in\mathtt{Trav}(Y) ] \ge
\mathbf{P}[t \in\mathtt{Trav}(X)]/(2D)$.
Now the rest of the proof goes through as when $L > 1$, with the small
change that instead of injectivity, we have that
$|\psi^{-1}(n)| \le2$ for $n \in\mathtt{Trav}(Y)$.
\end{pf}


\section{Proofs of Theorems \texorpdfstring{\protect\ref{t.main}}{1.1}
and \texorpdfstring{\protect\ref{t.main2}}{1.4}}\label{s.proofmain}


Here, we remove from Theorem~\ref{t.AGV} the
upper bound on the degrees in $G$ that was assumed
and we weaken the assumption on
the nature of nontrivial cycle frequency.


Consider a finite path $\mathcal{P}= \langle{e_t; 1 \le t \le
n}\rangle$.
Say that a time $t$ is a
\textit{cycle time of} $\mathcal{P}$  if there exist $1 \le
s \le t \le u \le n$
such that $(e_s, e_{s+1}, \ldots, e_u)$ is a cycle.
If the cycle is required to be an NT-cycle, then
we will call $t$ an NT-cycle time, and likewise for other types of cycles.
Call a cycle \textit{fully nontrivial}  if it is a loop or is
nontrivial and
its first edge is not the reverse of its last edge.
Such cycles will be called FNT-cycles.
For a finite or infinite path $\mathcal{P}$, we denote by $\mathcal
{P}\mathord{\upharpoonright}n$ its
initial segment of $n$ edges.

We state a slightly different version of Theorems \ref{t.main} and \ref{t.main2} here. At the end of the section,
we shall deduce the theorems as originally
stated in Section~\ref{s.intro}.

\begin{theorem}\label{t.AGV2}
Suppose that $G$ is a graph all of whose degrees are at least 3.
Let $X = \langle{X_t; t \ge1}\rangle$ be simple random walk on $G$.
If with positive probability
the limsup frequency of NT-cycle times of $X \mathord{\upharpoonright}n$
is positive as $n \to\infty$ (i.e., the expected limsup frequency is
positive), then the same is
true for nonbacktracking random walk. If  $G$ is also
$d$-regular, then $\rho(G) > 2\sqrt{d-1}/d$.
\end{theorem}


For $n \in\mathbb{Z}^+$, $\alpha> 0$, and a path $\mathcal{P}$ of
length $> n$,
let $C_n(\alpha, \mathcal{P})$ be the indicator that the number of
NT-cycle times
of  $\mathcal{P}\mathord{\upharpoonright}n$ is ${>}\alpha n$.
We shall prove the following finitistic version of Theorem~\ref{t.AGV2}, which will
be useful to us later.

\begin{theorem}\label{t.AGVfinite}
There exist $\zeta, \gamma> 0$ with the following property.
Suppose that $G$ is a graph all of whose degrees are at least 3.
Then for all $n$ and $\alpha$,
\[
\mathbf{E} \bigl[C_n(\alpha, X) \bigl(1 - C_n \bigl(\hat
\alpha, \mathtt{NB}(X) \bigr) \bigr) \bigr] < 3 e^{-\zeta n},
\]
where
\[
\hat\alpha := \gamma\alpha/\log^2 (10{,}368/\alpha).
\]
There exists $\zeta' > 0$ such that if  $G$ is also $d$-regular,
$\aaa (G) > 1$, and
\[
\mathbf{E} \bigl[C_n(\alpha, X) \bigr] > \frac{c_o n e^{-\zeta' n}}{\aaa (G) - 1},
\]
where $c_o$ is as in Proposition~\ref{p.cogrlimit},
then $\rho(G) > 2\sqrt{d-1}/d$.
If  $G$ is $d$-regular and
\[
\limsup_{n \to\infty} \bigl[\mathbf{E} {C_n(\alpha, X)}
\bigr]^{1/n} = 1,
\]
then
\[
\rho(G) > \frac{\sqrt{d-1}}{d} \bigl((d-1)^{\hat\alpha/24} +
(d-1)^{-\hat
\alpha
/24} \bigr).
\]
\end{theorem}

%

We shall use the following obvious fact.

\begin{lemma}\label{l.concat-NB}
If $(e_1, \ldots, e_k)$ and $(f_1, \ldots, f_m)$ are paths without
backtracking, the head of $e_k$ equals the tail of $f_1$, and $e_k$ is not
the reverse of $f_1$, then $(e_1, \ldots, e_k, f_1, \ldots, f_m)$ is
a path
without backtracking.
\end{lemma}

We shall apply the following well-known lemma to intervals with integer
endpoints.

\begin{lemma}[(Vitali covering)]\label{l.vitali}
Let $I$ be a finite collection of
subintervals of~$\mathbb{R}$. Write $\|I\|$ for the sum of the lengths
of the intervals in $I$. Then there exists a subcollection $J$ of
$I$ consisting of pairwise disjoint intervals such that $\|J\| \ge
\|I\|/3$.
\end{lemma}

This lemma is immediate from choosing iteratively the largest interval
disjoint from previously chosen intervals.

The following is a simple modification of a standard bound on large
deviations.

\begin{lemma}\label{l.exp-tails}
Suppose that $c > 0$.
There exist $\varepsilon\in(0, 1)$ and $\beta> 0$ such that
whenever $Z_1, \ldots, Z_n$ are random variables
satisfying the inequalities
$P[Z_k > z \mid Z_1, \ldots, Z_{k-1}] \le e^{-c z}$ for all $z > 0$,
we have
\[
\mathbf{P} \Biggl[\,\sum_{k=1}^{\varepsilon n}
Z_k \ge n \Biggr] \le e^{-\beta n}.
\]
\end{lemma}

\begin{pf}
Write $S_t := \sum_{j=1}^{\lfloor t\rfloor} Z_j$.
Given $\lambda:= c/2$ and $k \in[1, n]$, we have
\begin{eqnarray*}
\mathbf{E} \bigl[e^{\lambda Z_k} | Z_1, \ldots,
Z_{k-1} \bigr] &=& \int_0^\infty\mathbf{P}
\bigl[e^{\lambda Z_k} > z | Z_1, \ldots, Z_{k-1}
\bigr] \,dz \\
&\le &  1 + \int_1^\infty z^{-c/\lambda} \,dz
= 2,
\end{eqnarray*}
whence
\[
\mathbf{E} \bigl[e^{\lambda S_k} | Z_1, \ldots,
Z_{k-1} \bigr] \le 2 e^{\lambda S_{k-1}}.
\]
By induction, therefore, we have that $\mathbf{E} [e^{\lambda S_k} ]
\le2^k$.
It follows by Markov's inequality that
\[
\mathbf{P} \Biggl[\,\sum_{k=1}^{\varepsilon n}
Z_k \ge n \Biggr] = \mathbf{P} \bigl[e^{\lambda S_{\varepsilon n}} \ge
e^{\lambda n} \bigr] \le 2^{\varepsilon n} e^{-\lambda n}.
\]
Thus, if we choose $\varepsilon:= \min\{1/4, c/(4\log2)\}$ and $\beta:=
c/4$, the desired bound holds.
\end{pf}

Several lemmas now follow that will be used to handle various possible
behaviors of simple random walk on $G$.

\begin{lemma}\label{l.escape-tails}
Suppose that $G$ is a graph all of whose degrees are at least 3.
Let $X$ record the oriented edges taken by simple random walk on $G$.
Let $\Phi(n)$ index the $n$th edge of $X$ that remains in $\mathtt{NB}(X)$, so that
$\mathtt{NB}(X) = \langle{X_{\Phi(n)}}\rangle$: see~(\ref{e.defPhi}). Write $\Phi(0) := 0$.
Then there exists
$t_0 < \infty$ such that for all $n$ and all $t > t_0$, we have
\[
\mathbf{P} \bigl[\Phi(n+1) - \Phi(n) > t \bigr] < (8/9)^{t/2}.
\]
In addition, there exists $r$ such
that for every $n$ and $\lambda$,
\[
\mathbf{P} \bigl[\Phi(n) > (r+ \lambda) n \bigr] < (8/9)^{\lambda n/4}.
\]
More generally, for all $L > 0$,
there exists $r_L \le36^2 (8/9)^{L/4}$ such that
\[
\mathbf{P} \biggl[\,\sum_{k < n} \bigl(\Phi(k+1)-
\Phi(k) \bigr) \mathbf{1}_{[\Phi(k+1)-\Phi(k) >
L]} > (r_L + \lambda)n \biggr] <
(8/9)^{\lambda n/4}.
\]
\end{lemma}

\begin{pf}
Let $\widehat{X}$ be the lift of $X$ to $T$.
Then $\Phi$ also indexes the edges that remain in $\mathtt
{NB}(\widehat{X})$.
Since the distance of $\mathtt{NB}(\widehat{X})$ from $ \widehat{X}_1^-$ increases by 1
at each
step, the times $\Phi(n+1) - \Phi(n)$ are dominated by the times between
escapes for random walk on $\mathbb{N}$ that has probability $2/3$ to
move right and
$1/3$ to move left, reflecting at 0.
These in turn are dominated by the time to the first escape
for random walk~$S$ on $\mathbb{Z}$ with the same bias.
Such an escape can happen only at an odd time, $t$.
The chance of an escape at time exactly $t$ is
\begin{eqnarray*}
&&\mathbf{P}\bigl[S_{t-1} = 0, S_t = 1,
\forall t' > t \  S_{t'} > 0\bigr]
\\
&&\qquad= \pmatrix{t-1\cr (t-1)/2} \biggl(\frac{2}{3}\biggr)^{(t-1)/2}
\biggl(\frac{1}{3} \biggr)^{(t-1)/2}
\biggl(\frac{2}{3} \biggr) \biggl(\frac{1}{2} \biggr)
\\
&&\qquad \sim c \frac{(\sqrt{8/9} )^t}{\sqrt t}
\end{eqnarray*}
for some constant $c$.
This proves the first inequality.

Now $\Phi(n) = \sum_{k < n}  (\Phi(k+1) - \Phi(k) )$ and these
summands are dominated by the corresponding inter-escape times for the
biased random walk on $\mathbb{Z}$. The latter are i.i.d. with some distribution
$\nu$
(which we bounded in the last paragraph),
whence if we choose $a := (8/9)^{1/4} \in
 (\sqrt{8/9}, 1 )$ and put $b := \sum_{j \ge1} a^{-j} \nu
(j) \in
(1, \infty)$,
we obtain that for all $n$, we have $\mathbf{E} [a^{-\Phi(n)} ]\le b^n$.
By Markov's inequality, this implies that
\[
\mathbf{P} \bigl[\Phi(n) > (c + \lambda) n \bigr] \le a^{(c+\lambda)n}
b^n,
\]
so if we choose $c = r$ with $a^rb = 1$, then we obtain
the second inequality.

The third inequality follows similarly:
let $B_k :=  (\Phi(k+1) - \Phi(k) )\times\break
\mathbf{1}_{[\Phi(k+1) - \Phi(k) > L]}$.
Put $b_L := \nu[1,L] + \sum_{j > L} a^{-j} \nu(j) \in(1, 1 + 36 a^{L})$.
Then
$\mathbf{E} [a^{-\sum_{k=0}^{n-1} B_k} ]
\le b_L^n$
for all $n$.
By Markov's inequality, this implies that
\[
\mathbf{P} \Biggl[\,\sum_{k=0}^{n-1}
B_k > (c_L + \lambda) n \Biggr] \le a^{(c_L+\lambda)n}
b_L^n,
\]
so if we choose $c_L = r_L$ with $a^{r_L} b_L = 1$, then we obtain
the third inequality. We have the estimate $r_L \le36^2 a^{L}$.
\end{pf}

It follows that
\[
\mathbf{P} \bigl[\Phi\bigl(n/(r+\lambda)\bigr) > n \bigr] < a^{\lambda n}.
\]
That is, except for exponentially small probability, there are at least
$n/(r+\lambda)$ nonbacktracking edges by time $n$.
Similarly,
except for exponentially small probability, there are at most
$(r_L+\lambda)n$ edges by time $n$ that are in intervals of length${}>L$
that have no escapes.

The following is clear.

\begin{lemma}\label{l.NBin-cycle}
With notation as in Lemma~\ref{l.escape-tails}, if $s \le\Phi(n) \le
t$ satisfy
$ X_s^{-} =  X_t^{+}$, then $n$ is a cycle time of  $\mathtt{NB}(X)
\mathord{\upharpoonright}
t$.
\end{lemma}

We call\vspace*{1pt} a time $t$ an \textit{escape time}  for $\widehat{X}$ if
$-\widehat{X}_{t+1} \notin\mathtt{NB}(\widehat{X}_1, \ldots,
\widehat{X}_t)$ and $-\widehat{X}_{t+1}
\notin\{\widehat{X}_s; s > t+1\}$.
We let $\mathtt{Esc}(\widehat{X})$ be the set of escape times for
$\widehat{X}$.
Then $\mathtt{Esc}(\widehat{X}) = \operatorname{im}(\Phi-1)$.

\begin{lemma}\label{l.predictable}
Suppose that $\langle{\tau_k; 1 \le k < K}\rangle$ is a strictly increasing
sequence of stopping times for $X$, where $K$ is random, possibly
$\infty$.
Then there exist $\eta, \delta> 0$ such that for all $n > 0$,
\[
\mathbf{P} \bigl[K > n, \bigl|\bigl\{k \le n; \tau_k \in\mathtt {Esc}(
\widehat{X})\bigr\}\bigr| < \eta n \bigr] < e^{-\delta n}.
\]
\end{lemma}

\begin{pf}
Define random variables $\sigma_j$, $\lambda_j$ recursively.
First, we describe in words what they are.
Start by setting $\lambda_1 := 1$ and by
examining what happens after time $\tau_1$. If $\widehat{X}$ escapes,
define $\sigma_1 := 1$, $\lambda_2 := 2$,
and look at time $\tau_2$. If not, then
look at the first time $\tau_j$ that occurs after the first time $\ge t+1$
we know that $\widehat{X}$ has not escaped,\vspace*{1pt} that is, $t+1$ if
$-\widehat{X}_{t+1} \in\mathtt{NB}(\widehat{X}_1, \ldots, \widehat
{X}_t)$ or else
$\min\{ s > t+1; -\widehat{X}_{t+1} = \widehat{X}_s\}$,
and define $\sigma_1 := \tau_j - \tau_1$, $\lambda_2 := j$.
Now repeat from time $\tau_{\lambda_1}$ to define $\sigma_2$ and
$\lambda_3$, etc.

The precise definitions are as follows.
Suppose that $K > n$ (otherwise we do not define these random
variables).
Define $A_j$ to be the event that
one of the following holds:
\[
-\widehat{X}_{\tau_j+1} \in\mathtt{NB}(\widehat{X}_1, \ldots,
\widehat{X}_{\tau_j}) \quad\mbox{or}\quad \tau_j \in
\mathtt{Esc}(\widehat{X}).
\]
Write $\lambda_1 := 1$.
To recurse, suppose that $\lambda_k$ has been defined.
Let
\[
\lambda_{k+1} := \cases{ \lambda_k + 1, & \quad \mbox{if }$A_{\lambda_k}$,\cr
\min \bigl\{ j; -\widehat{X}_{\tau_{\lambda_k}+1}
\in\{\widehat{X}_s; \tau_{\lambda_k}+1 < s <
\tau_{j}\} \bigr\}, & \quad \mbox{otherwise}}
\]
and
\[
\sigma_k := \cases{ 1,  & \quad \mbox{if }$A_{\lambda_k}$,
\cr
\tau_{\lambda_{k+1}} - \tau_{\lambda_k},  &\quad \mbox{otherwise}.}
\]

Let $J := \max\{j; \lambda_j \le n\}$.
This is the number of times we have looked for escapes up to the
$n$th stopping time. Each stopping time until the $n$th is covered by one
of the intervals $[\tau_1, \tau_{\lambda_2}), \ldots, [\tau
_{\lambda_J},
\tau_{\lambda_{J+1}})$, which have lengths $\sigma_1, \ldots,
\sigma_J$.
Therefore,
we have that
\[
\sum_{j=1}^J \sigma_j \ge n.
\]
We claim that this forces $J$ to be large with high probability:
%
%
\begin{equation}\label{e.Jlarge}
\mathbf{P}[J \le\varepsilon n] \le e^{-\gamma n}
\end{equation}
for some $\varepsilon, \gamma> 0$.
Indeed, we claim that for each $k \le\varepsilon n$,
\[
\mathbf{P} \Biggl[\sum_{j=1}^k
\sigma_j \ge n \Biggr] \le e^{-\beta n},
\]
where $\varepsilon$ and $\beta$ are given by Lemma~\ref{l.exp-tails}
with $c$ (in that lemma) to be determined.
This would imply that
\[
\mathbf{P}[J \le\varepsilon n] \le \varepsilon n e^{-\beta n}.
\]

Now $\tau_{\lambda_{k+1}} - \tau_{\lambda_{k}} \ge\lambda_{k+1} -
\lambda_{k}$. Thus, it suffices to show that there is some $c > 0$ for
which
\[
\mathbf{P} [\lambda_{k+1} - \lambda_{k} \ge z |
\sigma_1, \ldots, \sigma_{k-1} ] \le e^{-c z}
\]
for all $z > 1$.
Now the event $\lambda_{k+1} - \lambda_{k} \ge z > 1$ implies the
event $B$
that $-\widehat{X}_t \notin\mathtt{NB}(\widehat{X}_1, \ldots,
\widehat{X}_{\tau_{\lambda_{k}}})$ for all $t \in(\tau_{\lambda_{k}},
\tau_{\lambda_{k}} + z)$ and that $-\widehat{X}_t =
\widehat{X}_{\tau_{\lambda_{k}}}$ for some $t \ge\tau_{\lambda
_{k}} + z$.
Because the distance from $\widehat{X}_t$ to $\widehat{X}_{\tau
_{\lambda_{k}}}$ has\vspace*{1pt} a
probability at least $2/3$ to get larger at all times, this is
exponentially unlikely in $z$.
What we need, however, is that this is exponentially unlikely even under
the given conditioning.
For every event $A$ in the $\sigma$-field on which we are
conditioning, we
always have that $A \supseteq[\tau_{\lambda_{k}} \in\mathtt{Esc}(\widehat{X})]$.
Furthermore, $\mathbf{P}[\tau_{\lambda_{k}} \in\mathtt
{Esc}(\widehat{X})] \ge1/2$. Hence, $\mathbf{P}(B
\mid A) \le2 \mathbf{P}(B)$, so that the bound on the unconditional
probability of
$B$ also gives an exponential bound on the conditional probability of $B$.
Thus, we have proved (\ref{e.Jlarge}).

Define $E_k := [\tau_{\lambda_k} \in\mathtt{Esc}(\widehat{X})]$.
We claim that
%
\begin{equation}\label{e.alwaysesc}
\mathbf{P} \bigl(E_k \mid\sigma(E_1, \ldots,
E_{k-1}) \bigr) \ge 1/2.
\end{equation}
Indeed, let $Z_t$ be the distance of $\widehat{X}_t^{+}$ to
$\widehat{X}_1^-$.
Note that $t \in\mathtt{Esc}(\widehat{X})$ iff $Z_s > Z_t$ for all
$s > t$.
Write $F_t(j)$ for the event that $Z_s > j$ for all $s > t$.
We claim that
\[
\mathbf{P} \bigl(E_k \mid\sigma(E_1, \ldots,
E_{k-1}, \widehat {X}_1, \ldots, \widehat{X}_{\lambda_k},
\lambda_1, \ldots, \lambda_k) \bigr) \ge 1/2,
\]
which is stronger than (\ref{e.alwaysesc}).
By choice of $\lambda_1, \ldots, \lambda_k$, we have that for every event
$E \in\sigma(E_1, \ldots, E_{k-1}, \widehat{X}_1, \ldots,
\widehat{X}_{\lambda_k}, \lambda_1, \ldots, \lambda_k)$,
\[
\mathbf{P}(E_k \mid E) = \mathbf{P} \bigl(F_{t}(j_m)
\mid F_{t}(j_1), \ldots, F_{t}(j_{m-1})
\bigr)
\]
for some $j_1, \ldots, j_{m-1} < j_m$ and some $t$, where $m \ge1$.
Since $F_t(j_i) \supseteq F_t(j_m)$ and $\mathbf{P} (F_t(j_m)
) \ge1/2$, the
claim follows.

Therefore, by (\ref{e.alwaysesc}),
we may couple the events $E_k$ to Bernoulli trials with
probability $1/2$ each so that the $k$th successful trial implies $E_k$. This
shows that there exists $\delta> 0$ such that
\[
\mathbf{P} \bigl[J > \varepsilon n, \bigl|\{k \le J; E_k\}\bigr| < \varepsilon
n/3 \bigr] < e^{-\delta n}.
\]
Hence,
\[
\mathbf{P} \bigl[K = \infty, \bigl|\bigl\{j \le n; \tau_j \in\mathtt
{Esc}(\widehat{X})\bigr\}\bigr| < \varepsilon n/3 \bigr] < e^{-\delta n}.
\]
Thus, we may choose $\eta:= \varepsilon/3$.
\end{pf}


Call a cycle of $>L$ edges an $L^+$-\textit{cycle}.
Define $I(n, L)$ to be the set of times $t \in[1, n]$ for which there
exist $1 \le s \le t \le u \le n$ such that $(X_s, X_{s+1}, \ldots, X_u)$
is a nontrivial $L^+$-cycle.

\begin{lemma}\label{l.whenlong}
If $n, L \ge1$ and
$\beta\in(0, 1)$, then
\[
\mathbf{E} \bigl[ \mathbf{1}_{[|I(n, L)| \ge\beta n]} \bigl(1 - C_n\bigl(
\beta/(2L), \mathtt{NB}(X)\bigr) \bigr) \bigr] < (8/9)^{\beta n/16}.
\]
\end{lemma}

\begin{pf}
Let $\mathtt{Long}$ be the event $ [|I(n, L)| \ge\beta n ]$.
Let $J$ be the union of intervals in $[1, n]$ that have length $> L$ and
are disjoint from $\mathtt{Esc}(X)$.
Let $\mathtt{Bad}$ be the event that $|J| > \beta n/2$.
By Lemma~\ref{l.escape-tails}, we have $\mathbf{P}(\mathtt{Bad}) <
(8/9)^{\beta n/16}$
(use $\lambda:= \beta/4$ there).
On the event $\mathtt{Long}\setminus\mathtt{Bad}$, the set $I(n, L)$ contains
at least $\beta n/2$ times that are within distance $L/2$ of an escape.
Therefore,
on the event $\mathtt{Long}\setminus\mathtt{Bad}$, there are at
least $\beta n/(2L)$
escapes in nontrivial cycles, whence $\mathtt{NB}(X)\mathord
{\upharpoonright}n$ has ${\ge}\beta n/(2L)$ NT-cycle times.
\end{pf}

Let $I_{\circ}(n, L)$ be the (random)
set of times $t \in[1, n] \setminus I(n, L)$
such that $X_t$ is a loop.\vspace*{-3pt}

\begin{lemma}\label{l.whenloop}
There exist $\eta, \delta> 0$ such that
if $n, L \ge1$ and
$\beta\in(0, 1)$, then
\[
\mathbf{E} \bigl[ \mathbf{1}_{[|I_{\circ}(n, L)| \ge\beta n]} \bigl(1 - C_n\bigl(\eta
\beta /(L+1), \mathtt{NB}(X)\bigr) \bigr) \bigr] < e^{-\delta n}.
\]
\end{lemma}

\begin{pf}
Let $\mathtt{Loop}:=  [|I_{\circ}(n, L)| \ge\beta n ]$.
Note that if there are 3 times at which a given loop in $G$ is
traversed, then necessarily the first of those times belongs to a
nontrivial cycle with at least one of the other times. In particular,
if a
given loop is traversed at least $L+2$ times, then it belongs
to a nontrivial long cycle. Therefore, on
$\mathtt{Loop}$, there are ${\ge}\beta n/(L+1)$ times spent in distinct loops.
If we take the first traversal of
a loop as a stopping time, then Lemma~\ref{l.predictable} supplies us
with $\eta,
\delta$ such that
on the event $\mathtt{Loop}$, except for probability $< e^{-\delta n}$,
the number of new loops at escape times is at least $\eta\beta n/(L+1)$.
Necessarily, all such loops remain in $\mathtt{NB}(X)$.
Therefore, $\mathtt{NB}(X)\mathord{\upharpoonright}n$ also has at
least $\eta\beta n/(L+1)$
loops
on the event $\mathtt{Loop}$ except for probability${}<e^{-\delta n}$.
\end{pf}

Let $D(n)$ denote the maximal number of disjoint FNT-cycles in $X
\mathord{\upharpoonright}n$, other than loops.\vspace*{-3pt}

\begin{lemma}\label{l.whenmany}
There exist $\eta, \delta> 0$ such that
if $n \ge1$ and
$\beta\in(0, 1)$, then
\[
\mathbf{E} \bigl[ \mathbf{1}_{[|D(n)| \ge\beta n]} \bigl(1 - C_n\bigl(\eta
\beta, \mathtt{NB}(X)\bigr) \bigr) \bigr] < e^{-\delta n}.
\]
\end{lemma}

\begin{pf}
Fix $n \ge1$.
Let $\mathtt{Cycs}$ be a (measurable) set of pairs of
times $1 \le s < t \le n$ such that $(X_s, X_{s+1}, \ldots, X_t)$
is an FNT-cycle other\vspace*{1pt} than a loop, chosen so that $|\mathtt{Cycs}| = D(n)$.
Let $\mathtt{Cycs\&Esc}:=  \{(s, t) \in\mathtt{Cycs}; t \in
\mathtt{Esc}(\widehat{X}) \}$.
By Lemma~\ref{l.predictable}, on the event $D(n) \ge\beta n$, we have
$|\mathtt{Cycs\&Esc}| > \eta'\beta n$ for some $\eta' > 0$
except for exponentially small probability.

Let $\mathtt{Sofar}$ be the set of $(s, t) \in\mathtt{Cycs}$ such
that $X_{s-1} \ne
- X_t$ or $s = 1$.

Note that for $(s, t) \in\mathtt{Cycs}$,
the cycle from $X_s$ to $X_t$ can be traversed in either
order, both being equally likely given $X_1, \ldots, X_{s-1}$, and at least
one of them has the property that $(s, t) \in\mathtt{Sofar}$ [see
Lemma~\ref{l.concat-NB}, where we concatenate $\mathtt{NB}(X_1, \ldots,
X_{s-1})$ with either
$\mathtt{NB}(X_s, \ldots, X_t)$ or $\mathtt{NB}(X_t, X_{t-1}, \ldots
, X_s)$, as
appropriate].
In fact, the same holds even conditioned on $\mathtt{Cycs\&Esc}$.
Therefore, we may couple to Bernoulli\vadjust{\goodbreak} trials and conclude that on the event
$D(n) \ge\beta n$, we have
$|\mathtt{Cycs\&Esc}\cap\mathtt{Sofar}| > \eta'\beta n/3$
except for exponentially small probability.
Note that for $t \in\mathtt{Cycs\&Esc}\cap\mathtt{Sofar}$, some
edge in the cycle
$(X_s, \ldots, X_t)$ belongs to $\mathtt{NB}(X)$ (see Lemma~\ref{l.concat-NB}
again)---more
precisely, $u \in\operatorname{im}\Phi$ for some $u \in[s, t]$---,
whence on the event $D(n) \ge\beta n$, we have
$\mathtt{NB}(X) \mathord{\upharpoonright}n$ has $> \eta\beta n$
cycle times
except for exponentially small probability, where $\eta:= \eta'/3$.
\end{pf}

\begin{lemma}\label{l.longcycles}
Suppose that $\rho:= \rho(G) < 1$, $n \ge2$, $\varepsilon> 0$ and $L
> 2e^2$.
Then
\[
\mathbf{P} \bigl[\bigl|\bigl\{\mbox{$L^+$-cycle times of } X
\mathord {\upharpoonright}n\bigr\}\bigr| \ge \varepsilon n \bigr] < e^{(6n/L)\log L}
\rho^{\varepsilon n/3}/(1 - \rho).
\]
\end{lemma}

\begin{pf}
For every $n$ and $k$,
the chance that $X_n$ begins a cycle of length $k$ is at most $\rho^k$.
Suppose that
the number of $L^+$-cycle times of
$X \mathord{\upharpoonright}n$ is at least $\varepsilon n$.
Then there are disjoint $L^+$-cycles in $X_1, \ldots, X_n$
the sum of whose lengths
is at least $\varepsilon n/3$ by Lemma~\ref{l.vitali}.
There are fewer than $n/L$ starting points and fewer than $n/L$ ending
points for those cycles since each has
length ${}>L$ and they are disjoint. The number of collections of
subsets of
$[0, n]$ of size at most $2n/L$ is ${<}e^{(n+1) h(2/L)} < e^{(6n/L)\log
L}$, where $h(\alpha) := -\alpha\log\alpha- (1-\alpha)\log
(1-\alpha)$.
This is because $h(\alpha) < -2\alpha\log\alpha$ for $\alpha< e^{-2}$.
For each such collection of starting and ending points giving
total length $k$, the chance that they do start
$L^+$-cycles is at most $\rho^k$, whence summing over collections and total
lengths that are ${\ge}\varepsilon n/3$, we get the result.
\end{pf}

Call a nonbacktracking cycle an \textit{NB-cycle}.
If an NB-cycle is a loop or
has the property that its last edge is different from the
reverse of its first edge, then call the cycle \textit{fully nonbacktracking},
abbreviated \textit{FNB-cycle}.
Recall that the number of NB-cycles of length $n$ starting from $x
\in\mathtt{V}(G)$ is $b_n(x)$.
We also say that a cycle starting from $x$ is ``at $x$''.
Let the number of FNB-cycles of length $n$ at $x$ be $b_n^*(x)$.
Recall that $S(x) := \{ n; b_n(x) \ne0\}$.
We shall need the following bounds on $b_n(x)$.

\renewcommand{\thetheorem}{1.5}
\begin{proposition}
Let $G$ be a graph with $\aaa (G) \ge1$.
For each $x \in\mathtt{V}(G)$, we have that
$\lim_{S(x) \ni n \to\infty} b_n(x)^{1/n}$ exists and
there is
a constant $c_x$ such that $b_n(x) \le c_x \aaa (G)^n$ for all
$n \ge1$.
Furthermore, if $x$ belongs to a simple cycle of length $L$, then $c_x
\le
2 + 2L\aaa (G)^{L-2}$.
If $G$ is $d$-regular, then $G$ is Ramanujan iff for all vertices $x$ and
all $n \ge1$, we have $b_n^*(x) \le2 (d-1)^{n/2}$.
\end{proposition}

\begin{pf}
Write $S^*(x) := \{ n; b_n^*(x) \ne0\}$.
Given two FNB-cycles starting at $x$,
we may concatenate the first with either
the second or the reversal of the second to obtain an FNB-cycle at $x$,
unless both FNB-cycles are the same loop.
Therefore, if $b_n^*(x)$ is the number of FNB-cycles at $x$, we have
$b_m^*(x) b_n^*(x)/2 \le b_{m+n}^*(x)$ for $m + n \ge3$, whence
$\langle{b_n^*(x)/2; n \ge2}\rangle$ is supermultiplicative and
Fekete's lemma
implies that $\lim_{S^*(x) \ni n \to\infty} b_n^*(x)^{1/n}$ exists and
$b_n^*(x) \le2 \aaa (G)^n$ for $n \ge2$.
It is easy to check that the same inequality holds for $n = 1$.
Together with Theorem~\ref{t.cogrowth},
this also implies that
if $G$ is $d$-regular and Ramanujan, then for all vertices $x$ and
all $n \ge1$, we have $b_n^*(x) \le2 (d-1)^{n/2}$.

Let $\widehat{b}_n(x) := b_n(x) - b_n^*(x)$ be the number of nonloop
NB-cycles at
$x$ whose last edge equals the reverse of its first edge, that is, NB-cycles
that are not FNB-cycles.
We shall bound $\widehat{b}_n(x)$ when
$x$ belongs to a simple cycle, say, $\mathcal{P}_0 = (e_1, \ldots, e_L)$
with $L$ edges.
Let $\mathcal{P}$ be a nonloop NB-cycle at $x$ whose last edge is
$e'$ and whose
first edge is $- e'$.
If $e'$ is a loop, then removing $e'$ at the end of $\mathcal{P}$
gives an FNB-cycle
$\mathcal{P}'$ at $x$.
Otherwise, decompose $\mathcal{P}$ as $\mathcal{P}_1.\mathcal{P}_2$, where
$.$ indicates concatenation, and $\mathcal{P}_2$ is maximal
containing only edges $e$ such
that $e \in\mathcal{P}_0$ or $- e \in\mathcal{P}_0$.
By reversing $\mathcal{P}_0$ if necessary, we may assume the former:
all edges of
$\mathcal{P}_2$ lie in~$\mathcal{P}_0$.
Suppose the first edge of $\mathcal{P}_2$ is $e_k$.
Then $\mathcal{P}_2$ traverses the remainder of $\mathcal{P}_0$ and
possibly the whole of
$\mathcal{P}_0$ several times.
Thus, write $\mathcal{P}_2 = \mathcal{P}_3.\mathcal{P}_4$, where
$\mathcal{P}_3 = (e_k, \ldots,
e_L)$ has length${}\le L$.
Finally, the NB-cycle $\mathcal{P}' := \mathcal{P}_1.(- e_{k-1},
\ldots,
- e_1).\overline{\mathcal{P}}_4$ is FNB, where the bar indicates
path reversal.
In addition, the length of $\mathcal{P}'$ differs from the length of
$\mathcal{P}$ by at
most $L-2$.
Since the map $\mathcal{P}\mapsto\mathcal{P}'$ is injective,
$\widehat{b}_n(x) \le\sum_{i=-1}^{L-2} b_{n+i}^*(x) \le2L \aaa (G)^{n+L-2}$.

Combining the results of the previous two paragraphs, we obtain that if $x$
belongs to a simple cycle, then there is a constant $c_x$ such that for all
$n \ge1$, we have $b_n(x) \le c_x \aaa (G)^n$.
We also get the bound claimed for $c_x$.

We now prove the same for
$x$ that do not belong to a simple cycle.
We claim that if $y$ is a neighbor of $x$, then $b_n(x) \le
b_{n-2}(y) + b_n(y) + b_{n+2}(y)$.
Indeed, let $\mathcal{P}$ be an NB-cycle at $x$.
Suppose the first edge of $\mathcal{P}$ goes to $y$.
If $\mathcal{P}$ is not FNB, then removing the first and last
edges of $\mathcal{P}$ yields an NB-cycle at $y$ of length $n-2$.
If $\mathcal{P}$ is FNB, then shifting the starting point from $x$ to $y$
yields an
FNB-cycle at $y$ of length $n$.
Lastly, if the first edge of $\mathcal{P}$ does not go to $y$, then we
may prepend
to $\mathcal{P}$ an edge from $y$ to $x$ and either append an edge
from $x$ to $y$
if the last edge of $\mathcal{P}$ was not from $y$, or else delete the last
edge of
$\mathcal{P}$, yielding an NB-cycle at $y$ of length $n+2$ or $n$.
This map of NB-cycles at $x$ to NB-cycles at $y$ is injective, which gives
the claimed inequality.
It follows that 
$b_n(x) \le c_x b_n(z)$, where
$z$ is the nearest point to $x$ that belongs to a simple cycle and $c_x$
does not depend on $n$.

Finally, if
$\lim_{S(x) \ni n \to\infty} b_n(x)^{1/n}$ exists for one $x$, then it
exists for all $x$ by the covering-tree argument we used earlier in
Section~\ref{s.intro}.
Suppose that for all $x$ belonging to a
simple cycle, $\lim_{S^*(x) \ni n \to\infty} b^*_n(x)^{1/n} <
\aaa (G)$.
Then the bounds in the preceding paragraphs show that
$\lim_{S^*(x) \ni n \to\infty} b_n(x)^{1/n} < \aaa (G)$.
It is not hard to see that therefore
$\limsup_{S(x) \ni n \to\infty} b_n(x)^{1/n} < \aaa (G)$ as
well, which
is a
contradiction to the definition of $\aaa (G)$.
Hence for some $x$, we have
$\lim_{S^*(x) \ni n \to\infty} b^*_n(x)^{1/n} = \aaa (G)$
and, therefore,
$\lim_{S(x) \ni n \to\infty} b_n(x)^{1/n} = \aaa (G)$ as well.
Together with Theorem~\ref{t.cogrowth},
this also implies that
if $G$ is $d$-regular and for all vertices $x$ and
all $n \ge1$, we have $b_n^*(x) \le2 (d-1)^{n/2}$, then $G$ is Ramanujan,
which completes the proof of
the last sentence of the proposition.
\end{pf}

Let $Y := \mathtt{NB}(X)$.
Let $A^L_n(\beta)$ be the event that there are $\ge\beta n$
times $t \in[1, n]$ for which there
exist $1 \le s \le t \le u \le n$ such that $(Y_s, Y_{s+1}, \ldots, Y_u)$
is a cycle with $u - s < L$.

\renewcommand{\thetheorem}{4.13}
\begin{lemma}\label{l.boost}
Let $G$ be $d$-regular with $\aaa (G) > 1$ and $\beta\in(0, 1)$.
For every $L < \infty$, if
\[
\mathbf{P} \bigl[A^L_n(\beta) \bigr] >
{c_o n (d-1)^{-\beta^2 n/6 + L/2} \over\aaa (G) - 1},
\]
where $c_o$ is as in Proposition~\ref{p.cogrlimit},
then $\rho(G) > 2\sqrt{d-1}/d$.
If
%
\begin{equation}\label{e.cond1}
\limsup_{n \to\infty} \mathbf{P} \bigl[A^{L}_n(
\beta) \bigr]^{1/n} = 1,
\end{equation}
then
\[
\rho(G) > {\sqrt{d-1} \over d} \bigl((d-1)^{\beta/24} +
(d-1)^{-\beta/24} \bigr).
\]
\end{lemma}

\begin{pf}
We may suppose that $\rho(G) < 1$, as there is nothing to prove
otherwise.

Let $A_n(\beta, L)$ be the event that $Y \mathord{\upharpoonright}n$ has
at least $\beta n$ cycle times
and that $Y_n$ completes a cycle of length ${\le}L$.
Note that $
\mathbf{P} [A_k(\beta, L) ]
\ge
\mathbf{P} [A^L_n(\beta) ]
/n$ for
some $k \in[\beta n, n]$ by considering the last cycle completed.

Consider the following transformation $\mathcal{P}\mapsto\mathcal
{P}'$ of finite
nonbacktracking paths $\mathcal{P}$:
let $I$ be the collection of cycles in $\mathcal{P}$. Choose
(measurably) a maximal
subcollection $J$ as
in Lemma~\ref{l.vitali}. Excise the edges in $J$ from $\mathcal{P}$,
concatenate the
remainder, and remove backtracks to arrive at $\mathcal{P}'$.
Then $\mathcal{P}'$ is a nonbacktracking path without cycles and
$|\mathcal{P}| - |\mathcal{P}'|$
is at least $1/3$ the number of cycle times of $\mathcal{P}$.


Fix $n$.
Let $p_n(\beta, L) := \mathbf{P} (A_n(\beta, L) )$.
Let $q_n(\beta, L)$ be the probability
that the length of $\mathcal{P}'$ is at most $n - \beta n/3$.
By the last paragraph, we have $q_n(\beta, L) \ge p_n(\beta, L)$.

We define another transformation $\mathcal{P}\mapsto\mathcal{P}''$
as follows, where
$\mathcal{P}''$ will be a nonbacktracking cycle when $Y_n$ completes
a cycle:
Let\vspace*{1.5pt} $m := \min\{i;  Y_i^{+} =  Y_n^{+}\}$.
Let $s$ be minimal with $\mathcal{P}'$
ending in $(Y_s, Y_{s+1}, \ldots, Y_n)$ and define $\widehat{\mathcal{P}}$ by
$\mathcal{P}' = \widehat{\mathcal{P}}.(Y_s, Y_{s+1}, \ldots, Y_n)$, where
$.$ indicates concatenation.
Since $\mathcal{P}'$ has no cycles, if $m < n$ (which it is if $Y_n$ completes
a cycle),
then $m < s$.
Now define $\mathcal{P}'' := \mathcal{P}.\overline{\mathcal{P}'}$
if $s = n$, where the bar indicates path reversal,
or else $\mathcal{P}'' := \mathcal{P}.(Y_{m+1}, Y_{m+2}, \ldots,
Y_s).\overline{\widehat{\mathcal{P}}}$.

Write $b := d-1$.
On the event
$A_n(\beta, L)$, we have that $\mathcal{P}''$ is a nonbacktracking cycle
with length at most $2n-\beta n/3 + L$. Furthermore, the map $\mathcal{P}
\mapsto\mathcal{P}''$ is injective because the first part of
$\mathcal{P}''$ is simply
$\mathcal{P}$.
Therefore, Proposition~\ref{p.cogrlimit} provides a constant $c_o$
such that
\[
d b^{n-1} q_n(\beta, L) \le \sum
_{k \le2n-\beta n/3 + L} c_o b_k(o),
\]
whence
\[
q_n(\beta, L) \le {c_o \aaa (G)^{2n-\beta n/3 + L} \over b^{n} (\aaa (G) - 1)}.
\]
%
For some $k \ge\beta n$,
we have
\[
{\mathbf{P} [A^L_n(\beta) ] \over n} \le p_k(\beta, L) \le q_k(
\beta, L) \le {c_o \aaa (G)^{2k-\beta k/3 + L} \over b^{k} (\aaa (G) - 1)}.
\]
It follows that if $\aaa (G) \le\sqrt b$, then
the last quantity above is
\[
\le {c_o b^{-\beta k/6 + L/2} \over\aaa (G) - 1} \le {c_o b^{-\beta^2 n/6 + L/2} \over\aaa (G) - 1},
\]
which proves the first part of the lemma.
Similarly, if
(\ref{e.cond1}) holds, then
\[
\aaa (G) \ge b^{1/(2-\beta/3)} > b^{1/2+\beta/12},
\]
whence by Theorem~\ref{t.cogrowth},
\[
\rho(G) > { b^{1/2+\beta/12} + b^{1/2-\beta/12} \over d}.
\]
\upqed\end{pf}

We remark that with more work, we may let $L := \infty$ in (\ref{e.cond1}).

\begin{pf*}{Proof of Theorem~\ref{t.AGVfinite}}
Let $X = \langle{X_n}\rangle$ be simple random walk on $G$ and
$\widehat{X}=
\langle{\widehat{X}_n}\rangle$ be its lift to the universal cover
$T$ of $G$.

Fix $n$.
Let $\mathtt{Good}$ be the event that $C_n(\alpha, X) = 1$.
We may choose $L \le34\log(10{,}368/\alpha)$ so that the number $r
_L$ of
Lemma~\ref{l.escape-tails} satisfies $r_L < \alpha/8$.
Fix such an $L$.


Let $\mathtt{Long}:=  [|I(n, L)| \ge\alpha n / 2 ]$.
By Lemma~\ref{l.whenlong} (using $\beta:= \alpha/2$), we have that
\[
\mathbf{E} \bigl[\mathbf{1}_{\mathtt{Long}} \bigl(1 - C_n\bigl(
\alpha /(4L), \mathtt{NB}(X)\bigr) \bigr) \bigr] < (8/9)^{\alpha n/32}.
\]


Let $\mathtt{Loop}:=  [|I_{\circ}(n, L)| \ge\alpha n/(8L) ]$.
Then by Lemma~\ref{l.whenloop} (using $\beta:= \alpha/4$),
\[
\mathbf{E} \bigl[{\mathbf{1}}_{\mathtt{Loop}} \bigl(1 - C_n\bigl(
\eta _1\alpha/\bigl(8L^2+8L\bigr), \mathtt{NB}(X)\bigr)
\bigr) \bigr] < e^{-\delta_1 n}
\]
for some $\eta_1, \delta_1 > 0$.

On the event $(\mathtt{Good}\setminus\mathtt{Long}\setminus\mathtt{Loop})$,
there are $\ge\alpha n/4$
times $t \in[1, n]$ for which there
exist $1 \le s \le t \le u \le n$ such that $(X_s, X_{s+1}, \ldots, X_u)$
is an NT-cycle with $1 \le u - s < L$ and that does not contain any loops;
this is because every loop can be contained in at most $2L$ NT-cycles of
length at most $L$ in $X \mathord{\upharpoonright}n$.
By Lemma~\ref{l.vitali},
on the event $(\mathtt{Good}\setminus\mathtt{Long}\setminus\mathtt{Loop})$,
there are ${\ge}\alpha n/(12 L)$ disjoint nonloop NT-cycles in $X
\mathord{\upharpoonright}
n$.
Within every NT-cycle, there is an FNT-cycle.
Thus, on the event $(\mathtt{Good}\setminus\mathtt{Long}\setminus
\mathtt{Loop})$,
there are ${\ge}\alpha n/(12 L)$ disjoint nonloop FNT-cycles in $X
\mathord{\upharpoonright}n$, that is, $D(n) > \alpha n/(12L)$
in the notation of Lemma~\ref{l.whenmany}.
Applying that lemma with $\beta:= \alpha/(12L)$ yields
\[
\mathbf{E} \bigl[{\mathbf{1}}_{\mathtt{Good}\setminus\mathtt
{Long}\setminus\mathtt{Loop}} \bigl(1 - C_n\bigl(
\eta_2 \alpha/(12L), \mathtt{NB}(X)\bigr) \bigr) \bigr] <
e^{-\delta_2 n}
\]
for some $\eta_2, \delta_2 > 0$.
Thus, the statement of the theorem holds with $\zeta:= \min
\{(\alpha/32)\log(9/8), \delta_1, \delta_2\}$ and $\gamma:= \min
\{\eta_1/12, \eta_2\}$.

%
%

Now we prove the second part of the theorem.

Suppose that $G$ is $d$-regular.
We may also suppose that $\rho(G) < (8/9)^{1/4}$, as there is nothing
to prove
otherwise.
Choose $L$ so that $L/\log L \ge2853/\alpha$, which is ${>}84/ (\alpha\log(1/\rho) )$.
Lemma~\ref{l.longcycles} then ensures that the above event $\mathtt
{Long}$ has
exponentially small probability:
\[
\mathbf{P}(\mathtt{Long}) < \frac{\rho^{\alpha n/84}}{1 - \rho}.
\]
%
Let $Y := \mathtt{NB}(X)$.
Although we did not state it, our proofs of Lemmas \ref{l.whenloop}
and~\ref{l.whenmany} provide many cycle times of $Y$ that occur in
cycles of length${}\le L$, that is, they show that the event
$A^L_n(\beta)$
occurs with high probability for certain $\beta$.
Thus,
\[
\mathbf{P} \bigl[\mathtt{Good}\setminus\mathtt{Long}\setminus
A^L_n(\hat\alpha) \bigr] < {\rho^{\alpha n/84} \over1 - \rho} +
e^{-\zeta n}.
\]
It follows by Lemma~\ref{l.boost} that if
\[
\mathbf{P}(\mathtt{Good}) \ge \frac{c_o n (d-1)^{-\alpha^2 n/24 + L/2}}{\aaa (G) - 1} + \frac{\rho^{\alpha n/84}}{1 - \rho}
+ e^{-\zeta n},
\]
then $\rho(G) > 2\sqrt{d-1}/d$.

Finally, if
$\limsup_{n \to\infty}  [\mathbf{E}C_n(\alpha, X) ]^{1/n}
= 1$, then
$\limsup_{n \to\infty}
\mathbf{P} [A^L_n(\hat\alpha,\break  Y) ]^{1/n} = 1$, so
Lemma~\ref{l.boost}
completes the proof.
\end{pf*}

\renewcommand{\thetheorem}{4.14}
\begin{remark}\label{r.rho}
Instead of requiring all degrees in $G$ to be at least 3, one could require
that $\rho(G) < 1$. A similar result holds.
\end{remark}

\begin{pf*}{Proof of Theorem~\ref{t.main2}}
Let $\mathcal{P}$ be an infinite path. Write $\alpha_n$ for the
number of NT-cycle
times${}\le n$ in $\mathcal{P}$, divided by $n$. Since we count here
cycles that may
end after time $n$, this may be larger than the density $\beta_n$
of NT-cycle times in $\mathcal{P}\mathord{\upharpoonright}n$.
However, we claim that $\limsup_{n \to\infty} \beta_n \ge\limsup_{n
\to\infty} \alpha_n$, whence the limsups are equal.

Suppose that $\alpha_n > \beta_n$.
Then there is some NT-cycle time $t \le n$ that belongs to an NT-cycle that
ends at some time $s > n$.
Every time in $[t, s]$ then is an NT-cycle time for $\mathcal{P}$.
It follows that $\beta_s \ge\alpha_n$, and this proves the claim.

It is now clear that Theorem~\ref{t.main2} follows from Theorem~\ref{t.AGV2}.
\end{pf*}

\begin{pf*}{Proof of Theorem~\ref{t.main}}
The proof follows just as for Theorem~\ref{t.main2}.
\end{pf*}

\section{Cycle encounters}\label{s.enc}

%

Here, we prove Theorem~\ref{t.prob10}.
We first sketch the proof that $q_n \to0$.
Assume the random walk has a good chance of encountering a short cycle
at a large time $n$. Because of the inherent fluctuations of
random walk, the time it reaches such a short cycle cannot be precise;
there must be many times around $n$ with approximately the same chance.
This means that there are actually many short cycles and if we look at
how many are encountered at times around $n$, we will have a good chance
of seeing many. This means the cycles are relatively dense (for random
walk) in that part of the graph, which boosts the cogrowth and hence the
spectral radius.


We begin by proving the following nonconcentration property of simple
random walk on regular graphs.\vspace*{-3pt}
\renewcommand{\thetheorem}{5.1}
\begin{lemma}\label{l.nonconcen}
Write $p_n(\cdot, \cdot)$ for the $n$-step transition probability of
simple random walk on a given graph.
Let $d < \infty$ and $\varepsilon> 0$. There exists $c > 0$ such that for
every $d$-regular graph $G$, every $o \in\mathtt{V}(G)$, and every $n
\ge1$, there exists $A \subseteq\mathtt{V}(G)$ that has the property\vspace*{-2pt} that
%
\begin{equation}\label{e.Abig}
p_n(o, A) > 1 - \varepsilon
\end{equation}
and\vspace*{-2pt}
%
\begin{equation}\label{e.Agood}
\forall x \in A,    \forall k \in [0, \sqrt n ]  \qquad p_{n+2k}(o, x) \ge
c p_n(o, x).
\end{equation}
\end{lemma}

\begin{pf}
Write $Q_n(j)$ for the probability that a binomial random variable with
parameters $\lfloor{n/2}\rfloor$ and $1/d$ takes the value $j$.
Given $\varepsilon$, define $c'$ so that\vspace*{-2pt}
\[
\sum_{|j - n/(2d)| \le c' \sqrt n} Q_n(j) > 1 -
\varepsilon^2.
\]
It has been known since the time of de Moivre\vspace*{-1pt} that
\[
Q_{n+2k}(j+k) \ge c Q_n(j)
\]
whenever $n \ge0$, $k \in[0, \sqrt n]$, and
$|j-n/(2d)| \le c' \sqrt n$.

Given $G$, $o\in\mathtt{V}(G)$, and $n \ge1$, let $X_1, \ldots,
X_n$ be
$n$ steps of simple random walk on $G$ starting with $ X_1^{-} = o$.
Define\vspace*{-2pt}
\[
Z := \bigl|\bigl\{i \in[1, n/2]; X_{2i-1} = - X_{2i}\bigr\}\bigr|.
\]
The events $ [X_{2i-1} = - X_{2i} ]$ are Bernoulli trials with
probability $1/d$ each, whence $Z$ has a binomial distribution with
parameters $\lfloor{n/2}\rfloor$ and $1/d$.
Thus,
\[
\mathbf{P} \bigl[\bigl|Z - n/(2d)\bigr| \le c' \sqrt n \bigr] > 1 -
\varepsilon^2
\]
by choice of $c'$.
Define\vspace*{-2pt}
\[
A := \bigl\{ x \in\mathtt{V}(G); \mathbf{P} \bigl[\bigl|Z - n/(2d)\bigr| \le
c' \sqrt n | X_{n}^{+} = x \bigr] > 1 -
\varepsilon \bigr\}.
\]
Since\vspace*{-6pt}
\begin{eqnarray*}
1 - \varepsilon^2 &<& \mathbf{P} \bigl[\bigl|Z - n/(2d)\bigr| \le
c' \sqrt n \bigr]
\\[-2pt]
&=& \sum_{ x \in\mathtt{V}(G)} p_n(o, x) \mathbf{P}
\bigl[\bigl|Z - n/(2d)\bigr| \le c' \sqrt n | X_{n}^{+}
= x \bigr]
\\[-2pt]
&\le& p_n(o, A) + \bigl(1 - p_n(o, A) \bigr) (1-
\varepsilon),
\end{eqnarray*}
we obtain (\ref{e.Abig}).

If we excise all even backtracking pairs $(X_{2i-1}, X_{2i})$ ($1 \le i
\le
n/2$) from the path $(X_1, \ldots, X_n)$, then we obtain simple random walk
for $n - 2Z$ steps conditioned not to have any even-time step be a backtrack.

Given $k \in[0, \sqrt n]$,
let $X'_1, \ldots, X'_{n+2k}$ be simple random walk from $o$
coupled with $X$ as follows:
Define
\[
Z' := \bigl|\bigl\{i \in\bigl[1, (n+2k)/2\bigr]; X'_{2i-1}
= - X'_{2i}\bigr\}\bigr|.
\]
By choice of $c$, we have $\mathbf{P}[Z' = j+k] \ge c \mathbf{P}[Z =
j]$ whenever $|j -
n/(2d)| \le c' \sqrt n$.
Thus, we may couple $X'$ and $X$ so that $Z' = Z+k$ with probability at
least~$c$ whenever $ X_n^{+} \in A$.
Furthermore, we may assume that the coupling is such that when $Z' = Z+k$
and we excise from each path the even backtracking pairs,
then what remains in $X'$ is the same as in $X$.
This implies that with probability at least $c$, we have ${X'}_{n+2k}^{+}
=  X_n^{+}$ whenever $X_n^+ \in A$.
This gives (\ref{e.Agood}).
\end{pf}

\renewcommand{\thetheorem}{1.2}
\begin{theorem}
Let $G$ be an infinite Ramanujan graph and $L \ge1$. Let
$q_n$ be
the probability that simple random walk at time $n$ lies on a nontrivial
cycle of length at most $L$. Then $\lim_{n \to\infty} q_n =
0$.
\end{theorem}

\begin{pf}
Let $S$ be the set of vertices that lie on a \textit{simple} cycle of
length at
most~$L$, so that $q'_n := \mathbf{P} [ X_n^- \in S ] =
\Omega(q_n)$ for
$n \ge~L$.
Suppose that $q'_n > 2\varepsilon$.
Choose $A$ and $c$ as in the lemma.
Then $\mathbf{P} [ X_n^- \in A \cap S ] \ge\varepsilon$, whence
$\mathbf{P} [ X_{n+2k}^- \in A \cap S ] \ge c \varepsilon$
for $k \in [0, \sqrt n ]$.

Let $I^L(n_1, n_2)$ be the number of
times $t \in[n_1, n_2]$ for which there
exist $n_1 \le s \le t \le u \le n_2$
such that $(X_s, X_{s+1}, \ldots, X_u)$
is a cycle with $u - s \le L$.
Then $\mathbf{E}_o [I^L(n, n + \sqrt n - 1),  X_n^- \in S
] \ge c'
\sqrt
n$ for some constant $c' > 0$ (depending only on $c \varepsilon$).
Thus, there is some vertex $x \in S$ (a value of $ X_{n}^-$) for which
$\mathbf{E}_x [I^L(1, \sqrt n) ] \ge c' \sqrt n$.
This also means
\[
\mathbf{P}_x \bigl[I^L(1, \sqrt n) \ge c'
\sqrt n/2 \bigr] \ge c'/2.
\]
Then Theorem~\ref{t.AGVfinite} completes the argument when $n$ is sufficiently
large since $c_x \le2+2L\aaa (G)^{L-2}$.
[The case $\aaa (G) = 1$ is immediate.]
Alternatively, one can appeal to Lemmas \ref{l.whenmany} and \ref{l.boost} instead of Theorem~\ref{t.AGVfinite}.
\end{pf}

\section*{Acknowledgement}
We are grateful to the referee for various
suggestions.

\printaddresses
\end{document}